\documentclass[11pt, twoside]{article}
\usepackage{cite}
\usepackage{amssymb}
\usepackage{mathrsfs}
\usepackage{amsmath}
\usepackage{amsthm}
\usepackage{amsfonts}
\usepackage{latexsym}
\usepackage{indentfirst}
\usepackage{color}
\usepackage{enumerate}
\usepackage[english]{babel}
\usepackage[colorlinks=true,
linkcolor=blue,
citecolor=red,
urlcolor=magenta]{hyperref}

\usepackage{txfonts}
\usepackage{anysize}

\allowdisplaybreaks

\newtheorem{theorem}{Theorem}[section]
\newtheorem{lemma}[theorem]{Lemma}
\newtheorem{corollary}[theorem]{Corollary}
\newtheorem{proposition}[theorem]{Proposition}

\theoremstyle{definition}

\numberwithin{equation}{section}

\begin{document}

\title{\bf\Large Sharp Weak-Type Inequalities
for Scalar and Matrix Weights
\footnotetext{\hspace{-0.35cm} 2020 {\it Mathematics Subject Classification}.
Primary 42B20; Secondary 47A30, 47A56, 42B35.\endgraf
{\it Key words and phrases.}
Calder\'on--Zygmund operator,
sparse domination,
weak-type inequality,
Muckenhoupt weight,
Cantor set.\endgraf
This project is partially supported by
the National Natural Science Foundation of China
(Grant Nos. 12431006 and 12371093),
the Beijing Natural Science Foundation (Grant No. 1262011), and
the Fundamental Research Funds for the Central Universities (Grant Nos.
2253200028 and 312200502406).}}
\date{}
\author{Fan Bu, Dachun Yang and Wen Yuan}

\maketitle

\vspace{-0.8cm}

\begin{center}
\begin{minipage}{13cm}
{\small {\bf Abstract}\quad
Let $p\in(1,\infty)$ and $p':=\frac p{p-1}$.
In this article, we establish a mixed multiplier
weak-type inequality with weight dependence
$[w]_{A_p}^{\frac1p}[w]_{A_\infty}^{\max\{\frac1p,\frac1{p'}\}}$
for scalar $A_p$ weights and Calder\'on--Zygmund operators
whose modulus of continuity satisfies the Dini condition.
As a consequence, we obtain the sharp exponent $\max\{1,\frac2p\}$
of $[w]_{A_p}$ in the corresponding pure $A_p$ bound.
Moreover, in the matrix-weighted setting,
we prove the sharpness of the exponent
$1+\min\{\frac1p,\frac1{p(p-1)}\}$ of $[W]_{\mathcal A_p}$
in the known multiplier weak-type inequality for the Hilbert transform.
}
\end{minipage}
\end{center}

\vspace{0.1cm}

\section{Introduction}

Muckenhoupt's characterization of the weighted
boundedness of the Hardy--Littlewood maximal operator \cite{m72}
is a foundational result in the theory of weighted norm inequalities.
Hunt et al. \cite{hmw73} subsequently established
the corresponding characterization for the Hilbert transform.
A function on $\mathbb R^n$ is called a \emph{scalar weight}
if it is positive almost everywhere and locally integrable.
For $p\in[1,\infty)$, a scalar weight $w$
belongs to the Muckenhoupt class $A_p$ if
\begin{equation}\label{eq:intro-ap}
[w]_{A_p}:=\begin{cases}
\displaystyle \sup_{\mathrm{cube}\,Q\subset\mathbb R^n}
\left(\fint_Qw\right)
\mathop{\operatorname{ess\,sup}}_{x\in Q} [w(x)]^{-1}
&\text{if } p=1, \\
\displaystyle \sup_{\mathrm{cube}\,Q\subset\mathbb R^n}
\left(\fint_Qw\right)\left(\fint_Q w^{-\frac1{p-1}} \right)^{p-1}
&\text{if } p\in(1,\infty)
\end{cases}
\end{equation}
is finite, where $\fint_Q w:=\frac1{|Q|}\int_Q w$.
Following earlier sharp estimates of Buckley \cite{b93} and
Petermichl \cite{p07,p08}, Hyt\"onen \cite{h12}
resolved the scalar $A_2$ conjecture by proving the linear
$A_2$ bound for general Calder\'on--Zygmund operators $T$.
This, together with sharp extrapolation \cite{dgpp05},
further implies that, for $p\in(1,\infty)$ and $w\in A_p$,
\begin{equation}\label{A2}
\|T\|_{L^p(w)\to L^p(w)}
\lesssim [w]_{A_p}^{\max\{1,\frac1{p-1}\}}.
\end{equation}
Lerner \cite{l13a} subsequently gave a simpler
using sparse domination approach \cite{l13b}.
For further developments in sparse domination,
see \cite{lor17,llo24,bd20,bbd22,bz26}.
The embedding $L^p\subset L^{p,\infty}$ yields a multiplier
weak-type estimate with the same power as the strong-type bound.
The question is whether this exponent can be improved in the weak-type setting.

Multiplier weak-type inequalities for the Hilbert transform and
the Hardy--Littlewood maximal operator were first studied by
Muckenhoupt and Wheeden \cite{mw77}.
Sawyer \cite{s85} later proved a mixed
weak-type inequality for the maximal operator on $\mathbb R$
involving two $A_1$ weights and conjectured an analogue result for
the Hilbert transform. Cruz-Uribe et al. \cite{cmp05}
resolved this conjecture and extended these inequalities to
the Hardy--Littlewood maximal operator and Calder\'on--Zygmund
operators on $\mathbb R^n$.
Let $T$ be a Calder\'on--Zygmund operator
and $M$ the Hardy--Littlewood maximal operator.
At the endpoint $p=1$, Cruz-Uribe et al.
\cite[Theorem~1.2]{cimpr21} proved that, for $w\in A_1$,
\begin{equation*}
\left\|wT\left( w^{-1}\cdot\right) \right\|_{L^1\to L^{1,\infty}}
\lesssim [w]_{A_1} [w]_{A_\infty}
\lesssim [w]_{A_1}^2,
\end{equation*}
where
$$
[w]_{A_\infty} := \sup_{\mathrm{cube}\,Q\subset\mathbb R^n}
\frac1{w(Q)}\int_QM(w\mathbf1_Q).
$$
Lerner et al. \cite[Theorem 1]{llor24b} showed that
the exponent $2$ is sharp for the Hilbert transform.
For $p\in(1,\infty)$, Cruz-Uribe and Sweeting
\cite[Theorem~1.5]{cs25} established the bound
\begin{equation*}
\left\|w^{\frac1p}T\left(w^{-\frac1p}\cdot\right)
\right\|_{L^p\to L^{p,\infty}}
\lesssim [w]_{A_p}^{\frac1p}[w]_{A_\infty}
\lesssim [w]_{A_p}^{1+\frac1p}.
\end{equation*}
For $p\in(1,2)$, Lerner et al. \cite[Theorem~1.2]{llor24a}
improved this bound to
\begin{equation}\label{two}
\left\|w^{\frac1p}T\left(w^{-\frac1p}\cdot\right)
\right\|_{L^p\to L^{p,\infty}}
\lesssim [w]_{A_p}^{1+\frac1{p^2}} [\log(e+[w]_{A_p})]^{\frac1p}.
\end{equation}
For $p\in[2,\infty)$, the $A_2$ theorem \eqref{A2}
and the embedding $L^p\subset L^{p,\infty}$ imply that
$$
\left\|w^{\frac1p}T\left(w^{-\frac1p}\cdot\right)
\right\|_{L^p\to L^{p,\infty}}
\lesssim [w]_{A_p}.
$$
Lerner et al. \cite[Theorem~1.4]{llor24a} proved that the exponent $1$
is sharp when $p\in[2,\infty)$.
For $p\in(1,2)$, they conjectured that \cite{llor24a}
the weight dependence in \eqref{two}
can be improved to $[w]_{A_p}^{\frac2p}$.

We now turn to matrix weights.
The theory of matrix-weighted inequalities originated
in the pioneering works \cite{tv97,nt96,v97}.
For subsequent developments concerning the boundedness of operators
in the matrix-weighted setting, see
\cite{n12,n13,nr18,n25,bc26,dly21,llsw26,yyz26}.
A \emph{matrix weight} is a matrix-valued function
$W:\ \mathbb R^n\to\mathbb C^{m\times m}$ with locally integrable entries
such that $W$ is positive definite almost everywhere.
or $p\in[1,\infty)$ and a measurable vector-valued function $\vec f$, define
\begin{equation*}
\|\vec f\|_{L^p(W)}
:=\left(\int_{\mathbb R^n} \left|W^{\frac1p}\vec f\right|^p\right)^{\frac1p}.
\end{equation*}
Following \cite{r03,fr04}, a matrix weight $W$ is called
an $\mathcal A_p$-matrix weight if
\begin{equation*}
[W]_{\mathcal A_p}:= \begin{cases}
\displaystyle \sup_{\mathrm{cube}\,Q\subset\mathbb R^n}
\mathop{\operatorname{ess\,sup}}_{y\in Q}
\fint_Q \left\|W^{\frac1p}(x)W^{-\frac1p}(y)\right\|^p\,dx
&\text{if } p=1, \\
\displaystyle \sup_{\mathrm{cube}\,Q\subset\mathbb R^n}
\fint_Q \left[ \fint_Q \left\|W^{\frac1p}(x)W^{-\frac1p}(y)\right\|^{p'}\,dy
\right]^{\frac{p}{p'}}\,dx
&\text{if } p\in(1,\infty)
\end{cases}
\end{equation*}
is finite. For $m=1$, this definition reduces to \eqref{eq:intro-ap}.
For $p\in(1,\infty)$, Cruz-Uribe et al. \cite[Corollary~1.16]{cim18} obtained
\begin{equation}\label{intro:matrix-strong-history}
\|T\|_{L^p(W)\to L^p(W)}
\lesssim [W]_{\mathcal A_p}^{1+\frac1{p(p-1)}}.
\end{equation}
The sharpness of the exponent $1+\frac1{p(p-1)}$
was proved in \cite{dptv24} for $p=2$
and in \cite{jqwx26,byy26} for all $p\in(1,\infty)$.

At the endpoint $p=1$, Cruz-Uribe et al. \cite[Theorem 1.2]{cimpr21} obtained
$$
\left\|WT\left(W^{-1}\cdot\right)
\right\|_{L^1\to L^{1,\infty}}
\lesssim [W]_{\mathcal A_1}^{2}.
$$
The exponent $2$ is sharp for the
Hilbert transform was established by Lerner et al. \cite[Theorem 1]{llor24b}.
For $p\in(1,\infty)$, Cruz-Uribe and Sweeting
\cite[Theorem~1.6]{cs25} obtained
$$
\left\|W^{\frac1p}T\left(W^{-\frac1p}\cdot\right)
\right\|_{L^p\to L^{p,\infty}}
\lesssim [W]_{\mathcal A_p}^{1+\frac1p}.
$$
This, together with \eqref{intro:matrix-strong-history}
and the embedding $L^p\subset L^{p,\infty}$, further implies that
\begin{equation}\label{cz:matrix-upper-main}
\left\|W^{\frac1p}T\left(W^{-\frac1p}\cdot\right)
\right\|_{L^p\to L^{p,\infty}}
\lesssim [W]_{\mathcal A_p}^{1+\min\{\frac1p,\frac1{p(p-1)}\}}.
\end{equation}
This leads to the question of whether the exponent in
\eqref{cz:matrix-upper-main} is sharp for all $p\in(1,\infty)$.

In this article, we prove the scalar conjecture of Lerner et al. \cite{llor24a}
and establish the sharpness of the matrix
exponent in \eqref{cz:matrix-upper-main} for the Hilbert transform.
To state the main results, we first recall some concepts.
A \emph{modulus of continuity} is a nondecreasing subadditive function
$\omega:[0,1]\to[0,\infty)$ with $\omega(0)=0$.
A bounded linear operator $T$ on $L^2(\mathbb R^n)$
is called an $\omega$-Calder'on--Zygmund operator
if it admits a kernel $K$, defined off the diagonal,
satisfying the following conditions:
\begin{enumerate}[\rm(i)]
\item for every compactly supported $f\in L^2$
and $x\notin\operatorname{supp}f$,
\begin{equation*}
Tf(x)=\int_{\mathbb R^n}K(x,y)f(y)\,dy,
\end{equation*}

\item there exists a positive constant $C_K$ such that,
for any $x,y\in\mathbb R^n$ with $x\ne y$,
\begin{equation*}
|K(x,y)|\le C_K |x-y|^{-n},
\end{equation*}

\item for any $x,y,h\in\mathbb R^n$ with $|x-y|>2|h|$,
\begin{equation*}
|K(x+h,y)-K(x,y)|+|K(x,y+h)-K(x,y)|
\le \omega\left(\frac{|h|}{|x-y|}\right) |x-y|^{-n}.
\end{equation*}
\end{enumerate}
We say that the modulus of continuity $\omega$
satisfies the \emph{Dini condition} if
\begin{equation*}
\|\omega\|_{\mathrm{Dini}}
:=\int_0^1\omega(t)\,\frac{dt}{t}<\infty.
\end{equation*}
Let
\begin{equation}\label{CT}
C_T:=\|T\|_{L^2\to L^2}+C_K+\|\omega\|_{\mathrm{Dini}}.
\end{equation}
The first main result in this article is the following mixed upper bound.

\begin{theorem} \label{thm:logfree-main}
Let $p\in(1,\infty)$ and $w\in A_p$.
Assume that the modulus of continuity $\omega$ satisfies the Dini condition
and $T$ is an $\omega$-Calder\'on--Zygmund operator.
Then there exists a positive constant $C$,
depending only on $n$ and $p$, such that
\begin{equation}\label{eq:logfree-main}
\left\|w^{\frac{1}{p}}T\left(w^{-\frac{1}{p}}\cdot\right)\right\|_{L^p\to L^{p,\infty}}
\le C C_T [w]_{A_p}^{\frac1p}[w]_{A_\infty}^{\max\{\frac1p,\frac1{p'}\}}.
\end{equation}
\end{theorem}

Since $[w]_{A_\infty}\lesssim[w]_{A_p}$ (see \cite[p.~778]{hp13}),
Theorem \ref{thm:logfree-main} yields
the following sharp bound in term of $[w]_{A_p}$ alone.

\begin{corollary} \label{cor:scalar-pure-sharp}
Let $p\in(1,\infty)$ and $w\in A_p$.
Assume that the modulus of continuity $\omega$ satisfies the Dini condition
and $T$ is an $\omega$-Calder\'on--Zygmund operator.
Then there exists a positive constant $C$,
depending only on $n$ and $p$, such that
\begin{equation*}
\left\|w^{\frac{1}{p}}T\left(w^{-\frac{1}{p}}\cdot\right)\right\|_{L^p\to L^{p,\infty}}
\le C C_T [w]_{A_p}^{\max\{1,\frac{2}{p}\}}.
\end{equation*}
\end{corollary}

Corollary \ref{cor:scalar-pure-sharp} answers the conjecture of Lerner et al. \cite{llor24a}.
Moreover, the sharpness of the exponent $\max\{1,\frac{2}{p}\}$ of $[w]_{A_p}$
in Corollary \ref{cor:scalar-pure-sharp} was established by
Lerner et al. \cite[Remark~3.3 and Theorem~1.4]{llor24a}.

The second main result in this article is the following theorem,
which proves the sharpness of the upper bound \eqref{cz:matrix-upper-main}.

\begin{theorem} \label{thm:hilbert-model}
Let $p\in(1,\infty)$ and $H$ be the Hilbert transform on $\mathbb R$.
Then there exists a positive constant $C$,
depending only on $p$, such that,
for any $m\ge2$ and $t\in[1,\infty)$,
\begin{equation*}
\sup_{[W]_{\mathcal A_p}\le t}
\left\|W^{\frac{1}{p}}H\left(W^{-\frac{1}{p}}\cdot\right)
\right\|_{L^p \to L^{p,\infty}}
\ge C t^{1+\min\{\frac1p,\frac1{p(p-1)}\}},
\end{equation*}
where the supremum is taken over all matrix weights
$W:\ \mathbb R\to\mathbb C^{m\times m}$ with
$[W]_{\mathcal A_p}\le t$.
\end{theorem}

The remainder of this article is organized as follows.

In Section \ref{sec:scalar}, we prove Theorem~\ref{thm:logfree-main}.
Following Lerner et al.
\cite[Proof of Theorem~1.2]{llor24a},
we reduce the problem to an estimate for $\|\mathcal A_{\mathcal S}
(w^{\frac1p}\mathbf1_E)\|_{L^{p'}(\sigma)}$
by using sparse domination and duality
(Corollary \ref{cor:sparse-testing-reduction}).
The equivalences of Cascante et al. \cite[Proposition~2.2]{cov04}
express this norm in terms of sums over descendant cubes
[\eqref{key} and \eqref{key2}].
A local estimate and the geometric decay of sparse generations
yield a density estimate with a logarithmic factor (Lemma~\ref{inner sum}).
The sharp reverse H\"older inequality supplies a positive power of the density
that absorbs this logarithm. This, together with dyadic Carleson embedding
and dyadic maximal estimates, further implies the desired bound
under the additional assumption $\Lambda_w(\mathcal S)\le2$
(Lemma \ref{prop:mixed-subfamily}).
Finally, using the pigeonholing technique,
we remove this assumption (Theorem~\ref{thm:layer}),
and hence complete the proof of Theorem~\ref{thm:logfree-main}.

In Section \ref{sec:matrix}, we prove Theorem \ref{thm:hilbert-model}
by adapting the matrix weight construction in
\cite{byy26} to the multiplier weak-type setting.
We modify the second scalar weight $w_1$
and the nondecreasing function $\varphi$ in that construction,
and choose a test function $\vec f$ whose weighted input norm
can be computed from a scalar integral.
For the output, a commutator with a nonnegative kernel
yields a lower bound on a test set of positive measure
at a fixed Cantor generation.
Positivity allows the contributions from pairwise disjoint
sibling intervals to add up on this set.
Comparing this lower bound with the input norm
gives the required operator norm estimate.

At the end of this introduction,
we make some notational conventions.
The notation $A\lesssim B$ means
that $A\leq CB$ for some positive constant $C$,
while $A\sim B$ means $A\lesssim B\lesssim A$.
We use $\mathbf 0$ to denote the origin of $\mathbb R^n$.
Let $\mathbb Z$ denote the set of all integers,
$\mathbb N:=\{1,2,\ldots\}$, and
$\mathbb Z_+:=\mathbb N\cup\{0\}$.
If $E$ is a measurable set in $\mathbb R^n$,
then we denote by $\mathbf 1_E$ its characteristic function.
For any $t\in\mathbb R$, we use the \emph{notation $\lceil t \rceil$}
to denote the smallest integer not less than $t$.
Finally, in all subsequent proofs we consistently retain the notation
introduced in the original theorem (or related statement).

\section{Upper Bounds for Scalar Weights}
\label{sec:scalar}

In this section, we prove the mixed scalar-weighted upper bounds
for $\omega$-Calder\'on--Zygmund operators
with $\|\omega\|_{\mathrm{Dini}}<\infty$
(Theorem~\ref{thm:logfree-main}).
As an application, we obtain the sharp scalar-weighted upper bound
(Corollary \ref{cor:scalar-pure-sharp}).
We begin with recalling some concepts.

A family $\mathcal D$ of cubes in $\mathbb R^n$ is called a
\emph{dyadic grid} if the following three conditions hold:
\begin{enumerate}[\rm(i)]
\item for every $Q\in\mathcal D$,
$\ell(Q)=2^{-k}$ for some $k\in\mathbb Z$;

\item for each $k\in\mathbb Z$, the collection
$\mathcal D_k:=\{Q\in\mathcal D:\ \ell(Q)=2^{-k}\}$
forms a partition of $\mathbb R^n$;

\item for any $Q,R\in\mathcal D$, either
$Q\cap R=\emptyset$, $Q\subset R$, or $R\subset Q$.
\end{enumerate}
A family $\mathcal S\subset\mathcal D$
is called \emph{$\frac12$-sparse} if
there exists a family of pairwise disjoint measurable sets
$\{E_Q\}_{Q\in\mathcal S}$ such that, for any $Q\in\mathcal S$,
$E_Q\subset Q$ and $|E_Q|\ge\frac12|Q|$.
For a $\frac12$-sparse family $\mathcal S\subset\mathcal D$,
its associated \emph{sparse operator} is defined by setting,
for any nonnegative measurable function $f$,
\begin{equation*}
\mathcal A_{\mathcal S}f
:=\sum_{Q\in\mathcal S}\left(\fint_Q f\right)\mathbf1_Q.
\end{equation*}

The following pointwise sparse domination
follows from \cite[Theorem~3.1 and Remark~4.3]{l16};
see also \cite[Theorem~4.2]{l17} and \cite[Theorem~2.4]{hrt17}.

\begin{lemma} \label{lem:sparse-domination}
Assume that the modulus of continuity $\omega$ satisfies the Dini condition
and $T$ is an $\omega$-Calder\'on--Zygmund operator.
Then there exists a positive constant $C$,
depending only on $n$, such that,
for any compactly supported $f\in L^1$,
there exist dyadic grids $\mathcal D_j$ and $\frac12$-sparse families
$\mathcal S_j\subset\mathcal D_j$, $j\in\{1,\ldots,3^n\}$, such that,
for almost every $x\in\mathbb R^n$,
\begin{equation*}
|Tf(x)|\le C C_T \sum_{j=1}^{3^n}\mathcal A_{\mathcal S_j}|f|(x),
\end{equation*}
where $C_T$ is as in \eqref{CT}.
\end{lemma}

The following equivalence follows from the duality
$(L^{p',1})^*=L^{p,\infty}$ and $(L^p)^*=L^{p'}$;
see \cite[Proof of Theorem~1.2]{llor24a}.

\begin{lemma} \label{lem:sparse-restricted-duality}
Let $p\in(1,\infty)$, $w\in A_p$, and $\sigma:=w^{-\frac1{p-1}}$.
Let $\mathcal S$ be a $\frac12$-sparse family
in a dyadic grid $\mathcal D$.
Suppose that there exists a positive constant $C$ such that,
for any measurable set $E$ with $|E|\in(0,\infty)$,
\begin{equation} \label{restricted testing estimate}
\left\|\mathcal A_{\mathcal S}
\left(w^{\frac1p}\mathbf1_E\right)\right\|_{L^{p'}(\sigma)}
\le C |E|^{\frac1{p'}}.
\end{equation}
Then, for any $f\in L^p$,
\begin{equation} \label{weak-type inequality}
\left\|w^{\frac1p}\mathcal A_{\mathcal S}
\left(w^{-\frac1p}|f|\right)\right\|_{L^{p,\infty}}
\le C\|f\|_{L^p}.
\end{equation}
Conversely, if $N$ is the best constant in \eqref{weak-type inequality},
then \eqref{restricted testing estimate} holds with $C$ replaced by $p'N$.
\end{lemma}

Combining Lemmas \ref{lem:sparse-domination} and
\ref{lem:sparse-restricted-duality},
we obtain the following corollary; we omit the details.

\begin{corollary}\label{cor:sparse-testing-reduction}
Let $p\in(1,\infty)$, $w\in A_p$, and $\sigma:=w^{-\frac1{p-1}}$.
Suppose that there exists a positive constant $C$, depending only on
$n$ and $p$, such that, for any dyadic grid $\mathcal D$,
any $\frac12$-sparse family $\mathcal S\subset\mathcal D$,
and any measurable set $E$ with $|E|\in(0,\infty)$,
\begin{equation}\label{eq:sparse-testing-target}
\left\|\mathcal A_{\mathcal S}
\left(w^{\frac1p}\mathbf1_E\right)\right\|_{L^{p'}(\sigma)}
\le C [w]_{A_p}^{\frac1p}[w]_{A_\infty}^{\max\{\frac1p,\frac1{p'}\}}|E|^{\frac1{p'}}.
\end{equation}
Then \eqref{eq:logfree-main} in Theorem~\ref{thm:logfree-main} holds
for every $\omega$-Calder\'on--Zygmund operator $T$
with $\|\omega\|_{\mathrm{Dini}}<\infty$.
\end{corollary}

Thus, proving Theorem~\ref{thm:logfree-main}
reduces to establishing the inequality \eqref{eq:sparse-testing-target}.
To show this inequality, following Lerner et al.
\cite[Proof of Theorem~1.2]{llor24a},
we apply the equivalence of Cascante et al. \cite[Proposition~2.2]{cov04},
which implies that, for any dyadic grid $\mathcal D$ and
any $\frac12$-sparse family $\mathcal S\subset\mathcal D$,
\begin{align}
\left\|\mathcal A_{\mathcal S}
\left(w^{\frac1p}\mathbf1_E\right)\right\|_{L^{p'}(\sigma)}^{p'}
&\sim \sum_{Q\in\mathcal S}
\left(\fint_Qw^{\frac1p}\mathbf1_E\right)\sigma(Q)
\left[\frac1{\sigma(Q)}\sum_{R\in\mathcal S(Q)}
\left(\fint_Rw^{\frac1p}\mathbf1_E\right)\sigma(R)
\right]^{p'-1} \label{key} \\
&\sim \left\| \sup_{Q\in\mathcal S} \mathbf 1_Q
\frac1{\sigma(Q)}\sum_{R\in\mathcal S(Q)}
\left(\fint_Rw^{\frac1p}\mathbf1_E\right)\sigma(R)
 \right\|_{L^{p'}(\sigma)}^{p'}, \label{key2}
\end{align}
where the positive equivalence constants depend only on $p$
and where $\mathcal S(Q):=\{R\in\mathcal S:\ R\subset Q\}$.
For $Q\in\mathcal S$, the generations of the sparse tree rooted at $Q$
are defined recursively. Let $\mathrm{ch}_{\mathcal S}^0(Q):=\{Q\}$.
For any $k\in\mathbb Z_+$,
$$
\mathrm{ch}_{\mathcal S}^{k+1}(Q)
:=\bigcup_{R\in \mathrm{ch}_{\mathcal S}^{k}(Q)}
\{P\in \mathcal S(R)\setminus\{R\}:\
P \text{ is a maximal element of } \mathcal S(R)\setminus\{R\}\}.
$$
Then $\mathcal S(Q)=\bigcup_{k=0}^{\infty}
\mathrm{ch}_{\mathcal S}^k(Q)$
and the cubes within each generation are pairwise disjoint.
To estimate the inner sum on the right-hand side of \eqref{key},
we need several technical lemmas.

\begin{lemma}\label{lem:mixed-packing}
Let $p\in(1,\infty)$, $w\in A_p$, and $\sigma:=w^{-\frac1{p-1}}$.
Suppose that $\mathcal S$ is a $\frac12$-sparse family
in a dyadic grid $\mathcal D$.
Then there exists a positive constant $C$, depending only on $p$,
such that, for any $Q\in\mathcal S$,
\begin{equation}\label{eq:mixed-packing}
\sum_{R\in\mathcal S(Q)}\left(\fint_Rw\right)\sigma(R)
\le C
\begin{cases}
\displaystyle [w]_{A_p(\mathcal S)}
\left(\fint_Q\sigma\right)^{2-p}|Q|
&\text{if } p\in(1,2],\\
\displaystyle [w]_{A_p(\mathcal S)}^{p'-1}
\left(\fint_Qw\right)^{2-p'}|Q|
&\text{if } p\in(2,\infty),
\end{cases}
\end{equation}
where
$$
[w]_{A_p(\mathcal S)}
:=\sup_{Q\in\mathcal S}
\left(\fint_Qw\right)\left(\fint_Q\sigma\right)^{p-1}.
$$
\end{lemma}

\begin{proof}
By \cite[Lemma~5.2]{h14}, we conclude that,
for any weight $v$, any $\theta\in[0,1)$, and $Q\in\mathcal S$,
\begin{equation}\label{eq:concave-packing}
\sum_{R\in\mathcal S(Q)}\left(\fint_Rv\right)^\theta|R|
\lesssim \left(\fint_Qv\right)^\theta|Q|.
\end{equation}
From the definition of $[w]_{A_p(\mathcal S)}$, it follows that,
for any $Q\in\mathcal S$ and $R\in\mathcal S(Q)$,
\begin{equation*}
\left(\fint_Rw\right)\sigma(R)\le
\begin{cases}
\displaystyle [w]_{A_p(\mathcal S)}
\left(\fint_R\sigma\right)^{2-p}|R|
&\text{if } p\in(1,2],\\
\displaystyle [w]_{A_p(\mathcal S)}^{p'-1}
\left(\fint_Rw\right)^{2-p'}|R|
&\text{if } p\in(2,\infty).
\end{cases}
\end{equation*}
Summing over $R\in\mathcal S(Q)$ and applying
\eqref{eq:concave-packing} with $(v,\theta)=(\sigma,2-p)$
in the first case and $(v,\theta)=(w,2-p')$ in the second,
we obtain \eqref{eq:mixed-packing}.
This completes the proof of Lemma \ref{lem:mixed-packing}.
\end{proof}

\begin{lemma}\label{lem:sparse-local-testing}
Let $p\in(1,\infty)$, $w\in A_p$, and $\sigma:=w^{-\frac1{p-1}}$.
Suppose that $\mathcal S$ is a $\frac12$-sparse family
in a dyadic grid $\mathcal D$.
Then there exists a positive constant $C$, depending only on $p$, such that
\begin{equation*}
\|\mathcal A_{\mathcal S}\sigma\|_{L^p(w)}^p
\le C [w]_{A_p(\mathcal S)} \sum_{Q\in\mathcal S}\sigma(Q).
\end{equation*}
\end{lemma}

\begin{proof}
Applying the equivalence \eqref{key}, we obtain
\begin{align*}
\|\mathcal A_{\mathcal S}\sigma\|_{L^p(w)}^p
&\sim \sum_{Q\in\mathcal S}
\left(\fint_Q\sigma\right)w(Q)
\left[ \frac1{w(Q)} \sum_{R\in\mathcal S(Q)}
\left(\fint_R\sigma\right)w(R) \right]^{p-1}\\
&= \sum_{Q\in\mathcal S} \left(\fint_Qw\right)\sigma(Q)
\left[ \frac1{w(Q)} \sum_{R\in\mathcal S(Q)}
\left(\fint_Rw\right)\sigma(R) \right]^{p-1}
=: \sum_{Q\in\mathcal S} I(Q).
\end{align*}
When $p\in(1,2]$, by Lemma~\ref{lem:mixed-packing}, we conclude that
\begin{align*}
I(Q)
&\lesssim \left(\fint_Qw\right)\sigma(Q)
\left[ \frac1{w(Q)} [w]_{A_p(\mathcal S)}
\left(\fint_Q\sigma\right)^{2-p}|Q| \right]^{p-1} \\
&= [w]_{A_p(\mathcal S)}^{p-1}
\left[\left(\fint_Qw\right)
\left(\fint_Q\sigma\right)^{p-1}\right]^{2-p}\sigma(Q)
\le [w]_{A_p(\mathcal S)}\sigma(Q).
\end{align*}
When $p\in(2,\infty)$, from Lemma~\ref{lem:mixed-packing} again, we deduce that
\begin{align*}
I(Q)
\lesssim \left(\fint_Qw\right)\sigma(Q)
\left[ \frac1{w(Q)} [w]_{A_p(\mathcal S)}^{p'-1}
\left(\fint_Qw\right)^{2-p'}|Q| \right]^{p-1}
= [w]_{A_p(\mathcal S)} \sigma(Q).
\end{align*}
This completes the proof of Lemma \ref{lem:sparse-local-testing}.
\end{proof}

\begin{lemma} \label{Carlson}
Let $\mathcal S\subset\mathcal D$ be a $\frac12$-sparse family in a
dyadic grid $\mathcal D$.
Then the following statements hold.
\begin{enumerate}[\rm(i)]
\item For any $Q\in\mathcal S$,
\begin{align*}
\sum_{R\in\mathcal S(Q)} |R| \le 2 |Q|.
\end{align*}

\item For any $Q\in\mathcal S$ and $k\in\mathbb Z_+$,
\begin{align*}
\sum_{R\in\mathrm{ch}_{\mathcal S}^k(Q)} |R|
\le \sum_{j=k}^{\infty} \sum_{R\in\mathrm{ch}_{\mathcal S}^j(Q)} |R|
\le 2^{1-k} |Q|.
\end{align*}
\end{enumerate}
\end{lemma}

\begin{proof}
Choose a family of pairwise disjoint measurable sets
$\{E_Q\}_{Q\in\mathcal S}$ such that, for any $Q\in\mathcal S$,
$E_Q\subset Q$ and $|E_Q|\ge\frac12|Q|$.
Then, for any $Q\in\mathcal D$,
\begin{align*}
\sum_{R\in\mathcal S(Q)} |R|
\le \sum_{R\in\mathcal S(Q)} 2|E_R|
\le 2 |Q|.
\end{align*}
This completes the proof of (i).

Now, we prove (ii).
For any $Q\in\mathcal S$ and $k\in\mathbb Z_+$, let
$$
T_k := \sum_{j=k}^{\infty} \sum_{R\in\mathrm{ch}_{\mathcal S}^j(Q)} |R|.
$$
By the just proven (i), we conclude that
\begin{align*}
T_k
=\sum_{R\in\mathrm{ch}_{\mathcal S}^k(Q)}
\sum_{P\in\mathcal S(R)}|P|
\le \sum_{R\in\mathrm{ch}_{\mathcal S}^k(Q)} 2|R|
=2(T_k-T_{k+1}),
\end{align*}
and hence $T_{k+1}\le\frac12T_k$.
This, together with (i) again, further implies that
\begin{align*}
T_k
\le 2^{-k} T_0
= 2^{-k} \sum_{R\in\mathcal S(Q)} |R|
\le 2^{1-k} |Q|,
\end{align*}
which completes the proof of (ii) and hence Lemma \ref{Carlson}.
\end{proof}

Now, we can estimate the inner sum on the right-hand side of \eqref{key}.

\begin{lemma} \label{inner sum}
Let $p\in(1,\infty)$, $w\in A_p$, and $\sigma:=w^{-\frac1{p-1}}$.
Let $\mathcal S\subset\mathcal D$ be a $\frac12$-sparse family in a
dyadic grid $\mathcal D$. Then there exists a positive constant $C$,
depending only on $p$, such that, for any measurable set $E$
and any $Q\in\mathcal S$ with $|E\cap Q|>0$,
\begin{equation}\label{eq:density-local-integral-logfree}
\sum_{R\in\mathcal S(Q)}\left(\fint_R w^{\frac1p}\mathbf1_E\right) \sigma(R)
\le C [w]_{A_p(\mathcal S)}^{\frac1p}
[\sigma(Q)]^{\frac1p}
|E\cap Q|^{\frac1{p'}}
\left(1+\log\frac{|Q|}{|E\cap Q|}\right)^{\frac1p}.
\end{equation}
\end{lemma}

\begin{proof}
Fix a measurable set $E$ and a cube $Q\in\mathcal S$ with $|E\cap Q|>0$.
Let $N\in \mathbb N$ be determined later, and let
$\mathcal H_{N,Q}:=\bigcup_{j=0}^{N-1}\mathrm{ch}_{\mathcal S}^j(Q)$.
By H\"older's inequality, we conclude that
\begin{align*}
I_N
:= \sum_{R\in\mathcal H_{N,Q}}
\left(\fint_Rw^{\frac1p}\mathbf1_E\right)\sigma(R)
= \int_{E\cap Q} w^{\frac1p} \mathcal A_{\mathcal H_{N,Q}}\sigma
\le |E\cap Q|^{\frac1{p'}}
\|\mathcal A_{\mathcal H_{N,Q}} \sigma\|_{L^p(w)}.
\end{align*}
This, together with Lemma~\ref{lem:sparse-local-testing}, further implies that
\begin{align}\label{IN}
I_N
\lesssim |E\cap Q|^{\frac1{p'}} [w]_{A_p(\mathcal S)}^{\frac1p}
\left[ \sum_{R\in\mathcal H_{N,Q}}\sigma(R) \right]^{\frac1p}
\le |E\cap Q|^{\frac1{p'}} [w]_{A_p(\mathcal S)}^{\frac1p}
N^{\frac1p} \sigma(Q)^{\frac1p} ,
\end{align}
where the last inequality used the fact that
$\sum_{R\in\mathcal H_{N,Q}}\mathbf1_R\le N\mathbf1_Q$.

For the remaining generations,
from H\"older's inequality and the definition of
$[w]_{A_p(\mathcal S)}$, we deduce that, for every $R\in\mathcal S(Q)$,
\begin{equation*}
\left(\fint_Rw^{\frac1p}\mathbf1_E\right)\sigma(R)
\le\left(\fint_Rw\right)^{\frac1p}\sigma(R)
\le [w]_{A_p(\mathcal S)}^{\frac1p}
[\sigma(R)]^{\frac1p} |R|^{\frac1{p'}}.
\end{equation*}
This, together with H\"older's inequality
and Lemma \ref{Carlson}(ii), further implies that
\begin{align}\label{IIN}
\sum_{j=N}^\infty \sum_{R\in\mathrm{ch}_{\mathcal S}^j(Q)}
\left(\fint_Rw^{\frac1p}\mathbf1_E\right)\sigma(R)
&\lesssim [w]_{A_p(\mathcal S)}^{\frac1p}
\sum_{j=N}^\infty \sum_{R\in\mathrm{ch}_{\mathcal S}^j(Q)}
[\sigma(R)]^{\frac1p} |R|^{\frac1{p'}} \notag \\
&\le [w]_{A_p(\mathcal S)}^{\frac1p} [\sigma(Q)]^{\frac1p}
\sum_{j=N}^\infty \left[ \sum_{R\in\mathrm{ch}_{\mathcal S}^j(Q)} |R|
\right]^{\frac1{p'}} \notag \\
&\le [w]_{A_p(\mathcal S)}^{\frac1p} [\sigma(Q)]^{\frac1p}
\sum_{j=N}^\infty (2^{1-j} |Q|)^{\frac1{p'}} \notag \\
&\sim [w]_{A_p(\mathcal S)}^{\frac1p} [\sigma(Q)]^{\frac1p}
(2^{1-N} |Q|)^{\frac1{p'}}.
\end{align}

Finally, choose
\begin{equation*}
N=\left\lceil\log_2\frac{2|Q|}{|E\cap Q|}\right\rceil.
\end{equation*}
Then $2^{1-N} |Q|\le|E\cap Q|$ and
\begin{equation*}
1\le N
\le 2+\log_2\frac{|Q|}{|E\cap Q|}
\le 2\left(1+\log\frac{|Q|}{|E\cap Q|}\right).
\end{equation*}
Using these, \eqref{IN}, and \eqref{IIN}, we obtain
\eqref{eq:density-local-integral-logfree}.
This completes the proof of Lemma \ref{inner sum}.
\end{proof}

To estimate $\|\mathcal A_{\mathcal S}(w^{\frac1p}\mathbf1_E)\|_{L^{p'}(\sigma)}$,
we still need several technical lemmas.
For any dyadic grid $\mathcal D$ and weight $w$,
the \emph{weighted dyadic maximal operator} $M_{\mathcal D,w}$ is
defined by setting, for any measurable function $f$ and $x\in\mathbb R^n$,
\begin{equation*}
M_{\mathcal D,w} f (x)
:=\sup_{\genfrac{}{}{0pt}{}{Q\in\mathcal D}{x\in Q}}
\frac1{w(Q)}\int_Q |f|w.
\end{equation*}
The following Hardy--Littlewood dyadic maximal theorem
is precisely \cite[Theorem~15.1]{ln19}.

\begin{lemma}
Let $w$ be a weight and $f$ a measurable function.
Then, for any $\lambda\in(0,\infty)$,
\begin{equation}\label{11}
\lambda w(\{x\in\mathbb R^n:\
M_{\mathcal D,w} f(x)>\lambda\})
\le \|f\|_{L^1(w)}
\end{equation}
and, for any $r\in(1,\infty]$,
\begin{equation}\label{pp}
\|M_{\mathcal D,w} f\|_{L^r(w)}
\le r' \|f\|_{L^r(w)}.
\end{equation}
\end{lemma}

The following dyadic Carleson embedding theorem
is a direct corollary of \cite[Theorem~4.5]{hp13}; we omit the details.

\begin{lemma} \label{embedding}
Let $r\in(1,\infty)$ and $w$ be a weight.
Suppose that $\mathcal S\subset\mathcal D$ is a $\frac12$-sparse family
in a dyadic grid $\mathcal D$. Then, for any $f\in L^r(w)$,
$$
\sum_{Q\in\mathcal S} w(Q) \left[\frac1{w(Q)}\int_Q|f|w\right]^r
\le(r')^r \Lambda_w(\mathcal S) \|f\|_{L^r(w)}^r,
$$
where $\Lambda_w(\mathcal S):=\sup_{R\in\mathcal S}
\frac1{w(R)}\sum_{Q\in{\mathcal S(R)}}w(Q)$.
\end{lemma}

We also need the following sharp reverse H\"older inequality,
which is precisely \cite[Theorem~2.3]{hpr12}.

\begin{lemma} \label{lem:mixed-rh}
Let $w\in A_\infty$
and $\delta:=\frac1{2^{n+1}[w]_{A_\infty}-1}$.
Then, for any cube $Q\subset\mathbb R^n$
and any $\varepsilon\in(0,\delta]$,
\begin{equation*}
\fint_Q w^{1+\varepsilon}
\le 2\left(\fint_Q w\right)^{1+\varepsilon}.
\end{equation*}
\end{lemma}

Now, we can estimate $\|\mathcal A_{\mathcal S}
(w^{\frac1p}\mathbf1_E)\|_{L^{p'}(\sigma)}$
under the assumption that $\Lambda_w(\mathcal S)\le2$.

\begin{lemma} \label{prop:mixed-subfamily}
Let $p\in(1,\infty)$, $w\in A_p$, and $\sigma:=w^{-\frac1{p-1}}$.
Suppose that $\mathcal S\subset\mathcal D$ is a $\frac12$-sparse family
in a dyadic grid $\mathcal D$ and $\Lambda_w(\mathcal S)\le2$.
Then there exists a positive constant $C$, depending only on $n$ and $p$,
such that, for every measurable set $E$ with $|E|\in(0,\infty)$,
\begin{equation*}
\left\|\mathcal A_{\mathcal S}
\left(w^{\frac1p}\mathbf1_E\right)\right\|_{L^{p'}(\sigma)}
\le C [w]_{A_p(\mathcal S)}^{\frac1p}
[w]_{A_\infty}^{\frac1p} |E|^{\frac1{p'}}.
\end{equation*}
\end{lemma}

\begin{proof}
Fix a measurable set $E$ with $|E|\in(0,\infty)$. Let
\begin{equation*}
\delta:=\frac1{2^{n+1}[w]_{A_\infty}-1}, \quad
r:=\frac{p'+1}{2},
\quad\text{and}\quad
\kappa:=\frac{\delta}{2(1+\delta)}
=\frac1{2^{n+2}[w]_{A_\infty}}.
\end{equation*}
By \eqref{key}, Lemma~\ref{inner sum},
and the definition of $[w]_{A_p(\mathcal S)}$, we obtain
\begin{align*}
\left\|\mathcal A_{\mathcal S}
\left(w^{\frac{1}{p}}\mathbf1_E\right)\right\|_{L^{p'}(\sigma)}^{p'}
\lesssim [w]_{A_p(\mathcal S)}^{\frac{p'}{p}}
\sum_{\genfrac{}{}{0pt}{}{Q\in\mathcal S}{|E\cap Q|>0}}
\frac{\int_{E\cap Q}w^{\frac{1}{p}}}{[w(Q)]^{\frac{1}{p}}}
|E\cap Q|^{\frac{1}{p}}
\left(1+\log\frac{|Q|}{|E\cap Q|}\right)^{\frac{p'-1}{p}}.
\end{align*}
This, together with H\"older's inequality, further implies that
\begin{align} \label{all est}
\left\|\mathcal A_{\mathcal S}
\left(w^{\frac{1}{p}}\mathbf1_E\right)\right\|_{L^{p'}(\sigma)}^{p'}
&\lesssim [w]_{A_p(\mathcal S)}^{\frac{p'}{p}}
\left\{ \sum_{\genfrac{}{}{0pt}{}{Q\in\mathcal S}{|E\cap Q|>0}}
\left[\frac{\int_{E\cap Q}w^{\frac{1}{p}}}{[w(Q)]^{\frac1{p}}}
\left(\frac{|E\cap Q|}{|Q|}\right)^{-\frac{\kappa}{p}}
\right]^{p'} \right\}^{\frac1{p'}} \notag\\
&\quad\times \left\{ \sum_{\genfrac{}{}{0pt}{}{Q\in\mathcal S}{|E\cap Q|>0}}
\left(\frac{|E\cap Q|}{|Q|}\right)^\kappa
|E\cap Q|
\left(1+\log\frac{|Q|}{|E\cap Q|}\right)^{p'-1}
\right\}^{\frac1p} \notag \\
&=: [w]_{A_p(\mathcal S)}^{\frac{p'}{p}}
\mathrm{I} \times \mathrm{II}.
\end{align}

From H\"older's inequality and the sharp reverse H\"older inequality
(Lemma \ref{lem:mixed-rh}), it follows that
\begin{equation}\label{intr}
w(E\cap Q)
\le \left( \int_Q w^{1+\delta} \right)^{\frac1{1+\delta}}
|E\cap Q|^{\frac{\delta}{1+\delta}}
\le 2w(Q) \left(\frac{|E\cap Q|}{|Q|}\right)^{\frac{\delta}{1+\delta}}
= 2w(Q) \left(\frac{|E\cap Q|}{|Q|}\right)^{2\kappa}.
\end{equation}
Using H\"older's inequality, \eqref{intr},
and the fact that $2\kappa(p'-r)=\frac{\kappa p'}{p}$, we conclude that
\begin{align*}
\left[\frac{\int_{E\cap Q}w^{\frac{1}{p}}}{[w(Q)]^{\frac1{p}}}
\left(\frac{|E\cap Q|}{|Q|}\right)^{-\frac{\kappa}{p}} \right]^{p'}
&\le\frac{(\int_{E\cap Q}w^{1-\frac{1}{r}})^r
[w(E\cap Q)]^{p'-r}}{[w(Q)]^{p'-1}}
\left(\frac{|E\cap Q|}{|Q|}\right)^{-\frac{\kappa p'}{p}}\\
&\le2^{p'-r}w(Q)
\left[\frac1{w(Q)}\int_Q
\left(w^{-\frac{1}{r}}\mathbf1_E\right)w\right]^r.
\end{align*}
This, together with the dyadic Carleson embedding theorem
(Lemma~\ref{embedding}) and $\Lambda_w(\mathcal S)\le2$,
further implies that
\begin{align}\label{I est}
\mathrm{I}^{p'}
\lesssim \Lambda_w(\mathcal S)
\left\|w^{-\frac{1}{r}}\mathbf1_E\right\|_{L^r(w)}^r
\lesssim |E|.
\end{align}

Next, we estimate $\mathrm{II}$.
By the change of variables $u:=\frac{\kappa}{2}\log\frac1t$
and $\kappa\in(0,1)$, we obtain
\begin{align*}
\sup_{t\in(0,1]} t^{\frac\kappa2}
\left( 1+\log \frac1t \right)^{p'-1}
&= \sup_{u\in[0,\infty)} e^{-u}
\left(1+\frac{2u}{\kappa}\right)^{p'-1}\\
&\le \left(\frac2\kappa\right)^{p'-1}
\sup_{u\in[0,\infty)} e^{-u}(1+u)^{p'-1}
\lesssim \kappa^{-(p'-1)},
\end{align*}
where the implicit positive constant depends only on $p$. Consequently,
\begin{align*}
\mathrm{II}^p
\lesssim \kappa^{-(p'-1)}
\sum_{Q\in\mathcal S}|Q|
\left(\frac{|E\cap Q|}{|Q|}\right)^{1+\frac\kappa2}.
\end{align*}
Choose a family of pairwise disjoint measurable sets
$\{E_Q\}_{Q\in\mathcal S}$ such that, for any $Q\in\mathcal S$,
$E_Q\subset Q$ and $|E_Q|\ge\frac12|Q|$.
From Cavalieri's principle and \eqref{11}, it follows that
\begin{align} \label{II est}
\mathrm{II}^p
&\lesssim \kappa^{-(p'-1)}
\sum_{Q\in\mathcal S}
\int_{E_Q}(M_{\mathcal D,1}\mathbf1_E)^{1+\frac\kappa2}
\le \kappa^{-(p'-1)}
\int_{\mathbb R^n}
(M_{\mathcal D,1}\mathbf1_E)^{1+\frac\kappa2} \notag \\
&= \kappa^{-(p'-1)} \left( 1+\frac\kappa2 \right)
\int_0^1 \lambda^{\frac\kappa2}
|\{M_{\mathcal D,1}\mathbf1_E>\lambda\}| \, d\lambda \notag \\
&\le \kappa^{-(p'-1)} \left( 1+\frac\kappa2 \right)
\int_0^1 \lambda^{\frac\kappa2-1} |E| \, d\lambda \notag \\
&= \kappa^{-(p'-1)} \left( 1+\frac\kappa2 \right) \frac2\kappa |E|
< 3 \kappa^{-p'} |E|,
\end{align}
where the last inequality used the fact that $\kappa\in(0,1)$.

Substituting \eqref{I est} and \eqref{II est} into
\eqref{all est}, we obtain
\begin{align*}
\left\|\mathcal A_{\mathcal S}
\left(w^{\frac{1}{p}}\mathbf1_E\right)\right\|_{L^{p'}(\sigma)}^{p'}
\lesssim [w]_{A_p(\mathcal S)}^{\frac{p'}{p}}
|E|^{\frac{1}{p'}}\left(\kappa^{-p'}|E|\right)^{\frac{1}{p}}
= [w]_{A_p(\mathcal S)}^{\frac{p'}{p}} \kappa^{-\frac{p'}{p}} |E|
\sim [w]_{A_p(\mathcal S)}^{\frac{p'}{p}} [w]_{A_\infty}^{\frac{p'}{p}} |E|,
\end{align*}
where the last equivalence used $\kappa^{-1}=2^{n+2}[w]_{A_\infty}$.
This completes the proof of Lemma \ref{prop:mixed-subfamily}.
\end{proof}

We now remove the additional assumption $\Lambda_w(\mathcal S)\le2$
by using the pigeonholing technique of \cite{lpr10}
(see also \cite{hptv14,h12,hp13}).

\begin{theorem} \label{thm:layer}
Let $p\in(1,\infty)$, $w\in A_p$, and $\sigma:=w^{-\frac1{p-1}}$.
Suppose that $\mathcal S\subset\mathcal D$ is a $\frac12$-sparse family
in a dyadic grid $\mathcal D$.
Then there exists a positive constant $C$, depending only on $n$ and $p$,
such that, for every measurable set $E$ with $|E|\in(0,\infty)$,
\begin{equation} \label{desired}
\left\|\mathcal A_{\mathcal S}
\left(w^{\frac1p}\mathbf1_E\right)\right\|_{L^{p'}(\sigma)}
\le C [w]_{A_p(\mathcal S)}^{\frac1p}
[w]_{A_\infty}^{\max\{\frac1p,\frac1{p'}\}} |E|^{\frac1{p'}}.
\end{equation}
\end{theorem}

\begin{proof}
Without loss of generality, we may assume that $\mathcal S$ is finite.
Decompose $\mathcal S$ into
\begin{equation}\label{eq:flat}
\mathcal S_k:=\left\{Q\in\mathcal S:\
2^k\le\left(\fint_Qw\right)\left(\fint_Q\sigma\right)^{p-1}
<2^{k+1}\right\},
\end{equation}
where $k\in\{0,\ldots,\lfloor\log_2[w]_{A_p(\mathcal S)}\rfloor\}$.
We claim that
\begin{equation}\label{claim}
\left\|\mathcal A_{{\mathcal S_k}}
\left(w^{\frac1p}\mathbf1_E\right)\right\|_{L^{p'}(\sigma)}
\lesssim 2^{\frac kp}[w]_{A_\infty}^{\max\{\frac1p,\frac1{p'}\}}|E|^{\frac1{p'}},
\end{equation}
where the implicit positive constant depends only on $n$ and $p$.
This, together with the triangle inequality and the geometric sum
\begin{equation*}
\sum_{k=0}^{\lfloor\log_2[w]_{A_p(\mathcal S)}\rfloor}2^{\frac kp}
\le\frac{[w]_{A_p(\mathcal S)}^{\frac1p}}{1-2^{-\frac1p}},
\end{equation*}
further implies the desired estimate \eqref{desired}.

Next, we prove \eqref{claim}.
Fix a measurable set $E$ with $|E|\in(0,\infty)$
and $k\in\{0,\ldots,\lfloor\log_2[w]_{A_p(\mathcal S)}\rfloor\}$.
We construct a principal family $\mathcal F\subset\mathcal S_k$.
Let $\mathcal F_0$ consist of the maximal cubes in $\mathcal S_k$ and,
for any $Q\in\mathcal S_k$, let
\begin{equation*}
\mathrm{ch}_{\mathcal F}(Q)
:=\operatorname{Max}_{\subset}
\left\{R\in\mathcal S_k(Q):\
\frac{1}{w(R)}\int_R w^{\frac1p}\mathbf1_E
>2\frac{1}{w(Q)} \int_Q w^{\frac1p}\mathbf1_E \right\},
\end{equation*}
where $\operatorname{Max}_{\subset}$ denotes the collection of maximal
elements with respect to inclusion.
For any $j\in\mathbb Z_+$, let
\begin{equation*}
\mathcal F_{j+1}
:=\bigcup_{Q\in\mathcal F_j}\mathrm{ch}_{\mathcal F}(Q)
\quad\text{and}\quad
\mathcal F:=\bigcup_{j=0}^\infty \mathcal F_j.
\end{equation*}
For any $Q\in\mathcal S_k$ with $|E\cap Q|>0$,
\begin{equation*}
\sum_{R\in\mathrm{ch}_{\mathcal F}(Q)}w(R)
\le\frac{w(Q)}{2\int_Q w^{\frac1p}\mathbf1_E}
\sum_{R\in\mathrm{ch}_{\mathcal F}(Q)}\int_R w^{\frac1p}\mathbf1_E
\le\frac12w(Q),
\end{equation*}
and hence
\begin{equation*}
\sum_{R\in\mathcal F(Q)}w(R)
\le \sum_{j=0}^\infty \frac1{2^j} w(Q)
= 2 w(Q).
\end{equation*}
Thus, $\Lambda_w(\mathcal F)\le2$.
By this, Lemma~\ref{prop:mixed-subfamily}, and \eqref{eq:flat}, we conclude that
\begin{equation}\label{AF}
\left\|\mathcal A_{\mathcal F}
\left(w^{\frac1p}\mathbf1_E\right)\right\|_{L^{p'}(\sigma)}
\lesssim
[w]_{A_p(\mathcal F)}^{\frac1p} [w]_{A_\infty}^{\frac1p} |E|^{\frac1{p'}}
\lesssim 2^{\frac kp}[w]_{A_\infty}^{\frac1p}|E|^{\frac1{p'}}.
\end{equation}

Now, we estimate $\|\mathcal A_{\mathcal S_k}
(w^{\frac1p}\mathbf1_E)\|_{L^{p'}(\sigma)}$.
The equivalence \eqref{key2} implies that
\begin{equation} \label{est}
\left\|\mathcal A_{\mathcal S_k}
\left(w^{\frac1p}\mathbf1_E\right)\right\|_{L^{p'}(\sigma)}
\sim \left\|\sup_{Q\in\mathcal S_k} \mathbf1_Q I(Q)
\right\|_{L^{p'}(\sigma)},
\end{equation}
where
\begin{equation*}
I(Q):=\frac1{\sigma(Q)} \sum_{R\in\mathcal S_k(Q)}
\left(\fint_R w^{\frac1p}\mathbf1_E\right) \sigma(R).
\end{equation*}
For any $R\in\mathcal S_k$, let $\pi(R)$ be the smallest cube in
$\mathcal F$ containing $R$. Then
$$
\frac{1}{w(R)}\int_R w^{\frac1p}\mathbf1_E
\le 2 \frac{1}{w(\pi(R))} \int_{\pi(R)} w^{\frac1p}\mathbf1_E.
$$
For any $R\in\mathcal S_k(Q)$, since $R\subset \pi(R)\cap Q$,
either $Q\subset \pi(R)$ or $\pi(R)\subsetneqq Q$, and hence
\begin{align} \label{est2}
\frac{1}{w(R)}\int_R w^{\frac1p}\mathbf1_E
\le 2 \inf_{x\in Q}
M_{\mathcal D,w} \left( w^{-\frac1{p'}}\mathbf1_E\right) (x)
+ 2 \sum_{\genfrac{}{}{0pt}{}{P\in\mathcal F}{R\subset P\subsetneqq Q}}
\frac{1}{w(P)} \int_P w^{\frac1p} \mathbf1_E.
\end{align}
By the definitions of $\mathcal S_k$ and $[w]_{A_p(\mathcal S_k)}$,
we conclude that, for any $Q\in\mathcal S_k$,
\begin{equation*}
[w]_{A_p(\mathcal S_k)}
\le 2^{k+1}
\le 2 \left(\fint_Qw\right)\left(\fint_Q\sigma\right)^{p-1}.
\end{equation*}
Substituting this estimate into \eqref{eq:mixed-packing},
we obtain, for any $Q\in \mathcal S_k$,
\begin{equation}\label{control}
\sum_{R\in\mathcal S_k(Q)}\left(\fint_Rw\right)\sigma(R)
\lesssim \left(\fint_Qw\right)\sigma(Q).
\end{equation}
From \eqref{est2}, \eqref{control}, and Tonelli’s theorem, it follows that
\begin{align*}
I(Q)
&\le \frac1{\sigma(Q)} \sum_{R\in\mathcal S_k(Q)} \left(\fint_R w\right)
\left[2 \inf_{x\in Q} M_{\mathcal D,w}
\left( w^{-\frac1{p'}}\mathbf1_E\right) (x)\right] \sigma(R) \\
&\quad+ \frac1{\sigma(Q)} \sum_{R\in\mathcal S_k(Q)} \left(\fint_R w\right)
\left[2 \sum_{\genfrac{}{}{0pt}{}{P\in\mathcal F}{R\subset P\subsetneqq Q}}
\frac{1}{w(P)} \int_P w^{\frac1p} \mathbf1_E\right] \sigma(R) \\
&\lesssim \left(\fint_Qw\right)\inf_{x\in Q}
M_{\mathcal D,w} \left( w^{-\frac1{p'}}\mathbf1_E\right) (x)
+\frac1{\sigma(Q)}
\sum_{\genfrac{}{}{0pt}{}{P\in\mathcal F}{P\subsetneq Q}}
\left(\fint_Pw^{\frac1p}\mathbf1_E\right)\sigma(P) \\
&\le \fint_Q M_{\mathcal D,w} \left( w^{-\frac1{p'}}\mathbf1_E\right) w
+ \frac1{\sigma(Q)}\int_Q
\mathcal A_{\mathcal F} \left( w^{\frac1p} \mathbf1_E \right)  \sigma.
\end{align*}
By this and \eqref{est}, we conclude that
\begin{align} \label{aim}
\left\|\mathcal A_{\mathcal S_k}
\left(w^{\frac1p}\mathbf1_E\right)\right\|_{L^{p'}(\sigma)}
&\lesssim \left\|\sup_{Q\in\mathcal S_k} \mathbf1_Q
\fint_Q M_{\mathcal D,w} \left( w^{-\frac1{p'}}\mathbf1_E\right)
w \right\|_{L^{p'}(\sigma)} \notag \\
&\quad+ \left\| M_{\mathcal D,\sigma} \left( \mathcal A_{\mathcal F} \left( w^{\frac1p} \mathbf1_E \right) \right) \right\|_{L^{p'}(\sigma)} \notag \\
&=:\mathrm{J}_1+\mathrm{J}_2.
\end{align}
For $\mathrm{J}_2$, \eqref{pp} and \eqref{AF} imply that
\begin{align*}
\mathrm{J}_2
\le p \left\| \mathcal A_{\mathcal F} \left( w^{\frac1p} \mathbf1_E \right) \right\|_{L^{p'}(\sigma)}
\lesssim 2^{\frac kp}[w]_{A_\infty}^{\frac1p}|E|^{\frac1{p'}}.
\end{align*}
For $\mathrm{J}_1$, choose a family of pairwise disjoint measurable sets
$\{E_Q\}_{Q\in\mathcal S}$ such that, for any $Q\in\mathcal S$,
$E_Q\subset Q$ and $|E_Q|\ge\frac12|Q|$.
Then, for any $Q\in \mathcal S_k$,
\begin{align*}
\sum_{R\in \mathcal S_k(Q)} w(R)
\le \sum_{R\in \mathcal S_k(Q)} 2 \int_{E_R} M(w\mathbf 1_Q)
\le 2 \int_Q M(w\mathbf 1_Q)
\le 2 [w]_{A_\infty} w(Q),
\end{align*}
and hence $\Lambda_w(\mathcal S_k)\le 2[w]_{A_\infty}$.
From this, \eqref{eq:flat}, and Lemma \ref{embedding}, it follows that
\begin{align*}
\mathrm{J}_1^{p'}
&\le \sum_{Q\in\mathcal S_k} \sigma(Q)
\left(\fint_Q w \right)^{p'} \left[\frac1{w(Q)} \int_Q
M_{\mathcal D,w} \left( w^{-\frac1{p'}}\mathbf1_E\right) w\right]^{p'} \\
&\le 2^{(k+1)(p'-1)}
\sum_{Q\in\mathcal S_k} w(Q)
\left[\frac1{w(Q)} \int_Q
M_{\mathcal D,w} \left( w^{-\frac1{p'}}\mathbf1_E\right) w\right]^{p'} \\
&\lesssim 2^{k(p'-1)}[w]_{A_\infty}
\left\|M_{\mathcal D,w} \left( w^{-\frac1{p'}}\mathbf1_E\right)
\right\|_{L^{p'}(w)}^{p'}
\lesssim 2^{k(p'-1)}[w]_{A_\infty}|E|.
\end{align*}
Substituting the estimates of $\mathrm{J}_1$ and $\mathrm{J}_2$
into \eqref{aim}, we obtain
\begin{equation*}
\left\|\mathcal A_{\mathcal S_k}\left(w^{\frac1p}\mathbf1_E\right)\right\|_{L^{p'}(\sigma)}
\lesssim2^{\frac kp}
\left([w]_{A_\infty}^{\frac1{p'}}+[w]_{A_\infty}^{\frac1p}\right)
|E|^{\frac1{p'}}
\lesssim2^{\frac kp}[w]_{A_\infty}^{\max\{\frac1p,\frac1{p'}\}}|E|^{\frac1{p'}}.
\end{equation*}
This proves \eqref{claim} and completes the proof of
Theorem~\ref{thm:layer}.
\end{proof}

Now, we can prove Theorem~\ref{thm:logfree-main}.

\begin{proof}[Proof of Theorem~\ref{thm:logfree-main}]
Combining Corollary \ref{cor:sparse-testing-reduction}
and Theorem \ref{thm:layer},
we obtain Theorem~\ref{thm:logfree-main}.
The factor $C_T$ enters only through
Lemma~\ref{lem:sparse-domination}; all remaining constants depend
only on $n$ and $p$.
\end{proof}

\section{Lower Bounds for Matrix Weights}
\label{sec:matrix}

In this section, we prove the sharp matrix-weighted
weak-type lower bounds for the Hilbert transform
(Theorem~\ref{thm:hilbert-model}).
We first introduce some notation that
will be used repeatedly throughout this section.

Starting from $\mathcal I_0:=\{[0,1]\}$, retain the two endpoint
subintervals of relative length $1/8$ from each interval in
$\mathcal I_j$ to form $\mathcal I_{j+1}$.
Thus $\mathcal I_j$ consists of $2^j$ intervals of length $8^{-j}$,
with total length $4^{-j}$ and mutual distances at least $6\cdot8^{-j}$.
Let
\begin{equation*}
E:=\bigcap_{j=0}^\infty\bigcup_{I\in\mathcal I_j}I, \quad
\rho(\cdot):=\operatorname{dist}(\cdot,E),
\end{equation*}
and, for any $j\in\mathbb Z_+$,
$$
\Omega_j:=\{x\in\mathbb R:\ 8^{-(j+1)}<\rho(x)\le8^{-j}\}.
$$
The Cantor set $E$ has measure zero.
For any $j\in\mathbb Z_+$ and
$I:=[b,b+8^{-j}]\in\mathcal I_j$, let
\begin{equation*}
H_I:=(b+3\cdot8^{-(j+1)},b+5\cdot8^{-(j+1)})
\quad\text{and}\quad
H_j:=\bigcup_{I\in\mathcal I_j}H_I.
\end{equation*}
Then $H_j\subset\Omega_j$ and $|H_j|=4^{-(j+1)}$.

Let $p\in(1,\infty)$,
\begin{equation*}
Q\in[\max\{10,10^{p-1}\},\infty),
\quad\text{and}\quad
\varepsilon:=Q^{-\frac1{p-1}}.
\end{equation*}
Then $\varepsilon\in(0,\frac1{10}]$.
For any $x\in\mathbb R$, let
\begin{equation} \label{phi}
\varphi(x):= \frac18 \int_{-\infty}^x
\sum_{j=0}^\infty 2^je^{\frac{\varepsilon}{p'}j}
\mathbf1_{\Omega_j\cap[0,1]}(t)\,dt.
\end{equation}
The $8^{-j}$-neighborhoods of the intervals in $\mathcal I_j$ cover $\Omega_j$,
which further implies that $|\Omega_j|\le3\cdot4^{-j}$, and hence
\begin{equation*}
\frac18\sum_{j=0}^\infty2^je^{\frac{\varepsilon}{p'}j}|\Omega_j\cap[0,1]|
\le\frac38\sum_{j=0}^\infty2^{-j}e^{\frac{\varepsilon}{p'}j}
=\frac38 \left( 1-\frac{e^{\varepsilon/p'}}2 \right)^{-1}
<1.
\end{equation*}
Thus $\varphi$ is absolutely continuous and nondecreasing,
$0\le\varphi<1$, and $\varphi\equiv0$ on $(-\infty,0]$.
Moreover, for almost every $x\in\mathbb R$,
\begin{align} \label{daoshu}
\varphi'(x)
=\frac18 \sum_{j=0}^\infty 2^je^{\frac{\varepsilon}{p'}j}
\mathbf1_{\Omega_j\cap[0,1]}(x).
\end{align}

For any $x\in\mathbb R$, let
\begin{equation*}
w_0(x):=
\begin{cases}
4^{-(p-1)j} e^{(p-1)\varepsilon j}
&\text{if } x\in\Omega_j \text{ for some } j\in\mathbb Z_+,\\
1&\text{if } \rho(x)>1 \text{ or } \rho(x)=0.
\end{cases}
\end{equation*}
For the second weight, let
\begin{equation*}
\delta:=Q^{-1}
\quad\text{and}\quad
N:=\lceil Q\rceil.
\end{equation*}
Then $\delta\in(0,\frac1{10}]$ and $\varepsilon^{p-1}=\delta$.
For any $x\in\mathbb R$, let
\begin{equation*}
w_1(x):=
\begin{cases}
\delta4^j &\text{if } x\in\Omega_j \text{ for some } j\in\{0,\ldots,N-1\},\\
4^N &\text{if } 0<\rho(x)\le8^{-N},\\
\delta &\text{if } \rho(x)>1 \text{ or } \rho(x)=0.
\end{cases}
\end{equation*}

The following construction of the matrix weight follows that in
\cite[Proposition 2.8]{byy26}, except that the definitions of
$w_1$ and $\varphi$ are modified to accommodate the weak-type estimate.

\begin{proposition}\label{direct:ap-bound}
Let $p\in(1,\infty)$,
\begin{equation*}
L:=\begin{pmatrix}
w_0^{\frac1p}&0\\
-\varphi w_1^{\frac1p}&w_1^{\frac1p}
\end{pmatrix},
\quad\text{and}\quad
W:=(L^*L)^{\frac p2}.
\end{equation*}
Then $W\in\mathcal A_p$ and
$[W]_{\mathcal A_p}\sim Q$,
where the positive equivalence constants depend only on $p$.
\end{proposition}

\begin{proof}
Fix a bounded interval $K$ of positive length. We claim that there exist
a generation $i\in\mathbb Z_+$ and a positive constant $C$, depending
only on $p$, such that
\begin{equation}\label{direct:claim-oscillation}
\operatorname{osc}_K\varphi
:= \sup_{x,y\in K} |\varphi(x)-\varphi(y)|
\le C 4^{-i} e^{\frac\varepsilon{p'}i},
\end{equation}
\begin{equation}\label{direct:claim-means}
\fint_K w_0 \le C 4^{-(p-1)i} e^{(p-1)\varepsilon i},
\quad
\fint_K w_0^{-\frac1{p-1}} \le C \varepsilon^{-1}4^ie^{-\varepsilon i},
\end{equation}
and
\begin{equation}\label{direct:modified-moments}
\fint_K w_1 \le C 4^{\min\{i,N\}},
\quad
\fint_K w_1^{-\frac1{p-1}} \le C
\delta^{-\frac1{p-1}}4^{-\frac{\min\{i,N\}}{p-1}}.
\end{equation}
Next, we prove this claim by considering three cases for $K$.

\textit{Case (1)} $|K|>\frac13$. In this case, take $i=0$.
The bounds $0\le\varphi<1$, $0<w_0\le1$, and $w_1\ge\delta$
imply \eqref{direct:claim-oscillation},
the first estimate in \eqref{direct:claim-means},
and the second estimate in \eqref{direct:modified-moments}.
Using $|K|>\frac13$, $|\Omega_j|\le3\cdot4^{-j}$,
and the definitions of $w_0$ and $w_1$, we obtain
\begin{align*}
\fint_Kw_0^{-\frac1{p-1}}
\le 1 + 3 \sum_{j=0}^\infty \int_{\Omega_j} w_0^{-\frac1{p-1}}
\le 1+9\sum_{j=0}^\infty e^{-\varepsilon j}
\sim 1 + \frac{1}{1-e^{-\varepsilon}}
\sim \varepsilon^{-1}
\end{align*}
and
\begin{align*}
\fint_K w_1
\le \delta + 3 \sum_{j=0}^\infty \int_{\Omega_j} w_1
\le\delta+9\delta N+9\cdot4^N\sum_{j=N}^\infty4^{-j}
\lesssim1,
\end{align*}
where the last inequality used $\delta N\le1+\delta$.

\textit{Case (2)} $|K|\le\frac13$ and
$\operatorname{dist}(K,E)\le2|K|$.
Choose the unique $i\in\mathbb Z_+$ such that
\begin{equation*}
8^{-(i+1)}<3|K|\le8^{-i}.
\end{equation*}
Then, for any $x\in K$,
$$
\rho(x)
\le\operatorname{dist}(K,E)+|K|
\le3|K|\le8^{-i}.
$$
Thus, $K\cap\Omega_j\neq\emptyset$ implies $j\ge i$.
For these $j$, enlarging $K$ by $8^{-j}$ at each endpoint gives
an interval of length at most $3\cdot8^{-i}$.
It meets at most one member of $\mathcal I_i$, and hence at most
$2^{j-i}$ generation $j$ intervals. Their $8^{-j}$-neighborhoods
cover $K\cap\Omega_j$. Thus, for any $j\in\mathbb Z_+$ with $j\ge i$,
\begin{equation}\label{direct:layer-cover}
|K\cap\Omega_j| \le 3\cdot2^{j-i}8^{-j}.
\end{equation}
Since $|K|>\frac1{24} 8^{-i}$, it follows that
\begin{equation}\label{direct:layer-proportion}
\frac{|K\cap\Omega_j|}{|K|}< 72\cdot4^{i-j}.
\end{equation}
By the definition of $\varphi$ and \eqref{direct:layer-cover},
\begin{align*}
\operatorname{osc}_K\varphi
&\le \frac18 \sum_{j=i}^\infty 2^j
e^{\frac{\varepsilon}{p'}j} |K\cap\Omega_j|
\le \frac38 \sum_{j=i}^\infty
e^{\frac{\varepsilon}{p'}j} 2^{-j-i} \\
&= \frac38 \frac{4^{-i}e^{\varepsilon i/p'}}{1-e^{\varepsilon/p'}/2}
< 4^{-i} e^{\frac\varepsilon{p'} i},
\end{align*}
which proves \eqref{direct:claim-oscillation}.
Since $4^{-(p-1)j}e^{(p-1)\varepsilon j}$ decreases with $j$, we have
\begin{align*}
\fint_K w_0 \le 4^{-(p-1)i} e^{(p-1)\varepsilon i}.
\end{align*}
Moreover, \eqref{direct:layer-proportion} implies that
\begin{equation*}
\fint_Kw_0^{-\frac1{p-1}}
= \frac1{|K|}\sum_{j=i}^\infty \int_{K\cap\Omega_j} w_0^{-\frac1{p-1}}
\le 72 \sum_{j=i}^\infty4^{i-j}4^je^{-\varepsilon j}
= 72 \frac{4^i e^{-\varepsilon i}}{1-e^{-\varepsilon}}
\sim \varepsilon^{-1}4^ie^{-\varepsilon i}.
\end{equation*}
These two estimates prove \eqref{direct:claim-means}.
For $w_1$, if $i\le N-1$, then \eqref{direct:layer-proportion} implies that
\begin{align*}
\fint_Kw_1
&= \frac1{|K|}\sum_{j=i}^\infty \int_{K\cap\Omega_j} w_1
\le 72 \delta\sum_{j=i}^{N-1}4^{i-j}4^j
+ 72\cdot4^N\sum_{j=N}^\infty4^{i-j} \\
&= 72\cdot4^i\left[\delta(N-i)+\frac43\right]
\lesssim 4^i,
\end{align*}
where the last inequality used $\delta N\le1+\delta$.
If $i\ge N$, then $w_1=4^N$ almost everywhere on $K$.
In both cases, $w_1\ge\delta4^{\min\{i,N\}}$ almost everywhere on $K$.
Therefore, \eqref{direct:modified-moments} holds.

\textit{Case (3)} $|K|\le\frac13$ and
$\operatorname{dist}(K,E)>2|K|$.
If $\operatorname{dist}(K,E)>\frac18$, take $i=0$.
Then all the claimed estimates
follow from $0\le\varphi<1$ and $w_0\equiv 1$ and $w_1\equiv\delta$ on $K$.
Otherwise, choose the unique $i\in\mathbb N$ such that
\begin{equation*}
8^{-(i+1)}<\operatorname{dist}(K,E)\le8^{-i}.
\end{equation*}
Then, for any $x\in K$,
\begin{equation*}
8^{-(i+1)}<\rho(x)
\le\operatorname{dist}(K,E)+|K|<\frac32\,8^{-i}.
\end{equation*}
Consequently, $K$ meets only $\Omega_{i-1}$ and $\Omega_i$,
and hence, for any $x\in K$,
\begin{align*}
4^{-(p-1)i}e^{(p-1)\varepsilon i}
\le w_0(x)
\le 4^{p-1} 4^{-(p-1)i} e^{(p-1)\varepsilon i}
\end{align*}
and
\begin{align*}
\frac\delta4 4^{\min\{i,N\}}
\le w_1(x)
\le 4^{\min\{i,N\}}.
\end{align*}
These pointwise estimates imply
\eqref{direct:claim-means} and \eqref{direct:modified-moments}.
Finally, from \eqref{daoshu}, it follows that,
for almost every $x\in K$,
$$
\varphi'(x) \le \frac18 2^i e^{\frac\varepsilon{p'} i}.
$$
This, together with the fact that $|K|<\frac12\operatorname{dist}(K,E) \le \frac12 8^{-i}$,
further implies that
\begin{equation*}
\operatorname{osc}_K\varphi
\le \frac18 2^i e^{\frac\varepsilon{p'} i} |K|
\le \frac1{16} 4^{-i} e^{\frac\varepsilon{p'} i},
\end{equation*}
which proves \eqref{direct:claim-oscillation}.
This finishes the proof of the claim.

Since $W^{\frac2p}=L^*L$, we have $|W^{\frac1p}\vec z|=|L\vec z|$
for all $\vec z\in\mathbb C^2$,
and hence, for any $x,y\in\mathbb R$,
\begin{equation*}
\left\|W^{\frac1p}(x)W^{-\frac1p}(y)\right\|
=\left\|L(x)L^{-1}(y)\right\|.
\end{equation*}
A direct calculation shows that
\begin{equation*}
L(x)L^{-1}(y)=
\begin{pmatrix}
[w_0(x)]^{\frac1p}[w_0(y)]^{-\frac1p}&0\\
[\varphi(y)-\varphi(x)]
[w_1(x)]^{\frac1p}[w_0(y)]^{-\frac1p}&
[w_1(x)]^{\frac1p}[w_1(y)]^{-\frac1p}
\end{pmatrix}.
\end{equation*}
From these, the equivalence of norms on finite-dimensional spaces,
and the just proven claim, we infer that
\begin{align*}
&\fint_K \left(\fint_K
\left\|W^{\frac{1}{p}}(x)W^{-\frac{1}{p}}(y)\right\|^{p'}\,dy
\right)^{p-1}\, dx \\
&\quad\lesssim \left(\fint_Kw_0\right)
\left(\fint_Kw_0^{-\frac1{p-1}}\right)^{p-1}
+\left(\fint_Kw_1\right)
\left(\fint_Kw_1^{-\frac1{p-1}}\right)^{p-1}\\
&\qquad+(\operatorname{osc}_K\varphi)^p
\left(\fint_Kw_1\right)
\left(\fint_Kw_0^{-\frac1{p-1}}\right)^{p-1}\\
&\quad\lesssim \varepsilon^{-(p-1)}+\delta^{-1}
+\varepsilon^{-(p-1)} 4^{\min\{i,N\}-i}
\lesssim Q,
\end{align*}
where the last inequality used $\varepsilon^{p-1}=\delta=Q^{-1}$.
Thus, $[W]_{\mathcal A_p}\lesssim Q$.

For the lower bound, take
$K:=(-\frac32\,8^{-N},-\frac12\,8^{-N})$ and denote its left
and right halves by $K_-$ and $K_+$.
Note that $w_1\equiv \delta4^{N-1}$ on $K_-$ and
$w_1\equiv 4^N$ on $K_+$. Thus,
$$
[W]_{\mathcal A_p}\ge
\left(\fint_K w_1\right)
\left(\fint_K w_1^{-\frac1{p-1}}\right)^{p-1}
> \left(\frac12 4^N \right)
\left[\frac12 (\delta4^{N-1})^{-\frac1{p-1}}\right]^{p-1}
=2^{2-p} \delta^{-1}=2^{2-p}Q.
$$
This completes the proof of Proposition \ref{direct:ap-bound}.
\end{proof}

Now, we can prove Theorem~\ref{thm:hilbert-model}.

\begin{proof}[Proof of Theorem~\ref{thm:hilbert-model}]
We first consider the case $m=2$. Note that
$$
\left\|W^{\frac{1}{p}}H\left(W^{-\frac{1}{p}}\cdot\right)
\right\|_{L^p \to L^{p,\infty}}
=\left\|W^{\frac{1}{p}}H \right\|_{L^p(W) \to L^{p,\infty}}.
$$
Let $W$ be as in Proposition \ref{direct:ap-bound},
\begin{equation*}
J:= N+\lceil\varepsilon^{-1}\rceil, \quad
u:= \varepsilon \sum_{j=0}^{J-1} 4^j
e^{-\frac{\varepsilon}{p'}j} \mathbf1_{H_j},
\quad\text{and}\quad
\vec f:= \begin{pmatrix}u\\
\varphi u\end{pmatrix},
\end{equation*}
where $\varphi$ is as in \eqref{phi}.
Both $u$ and $\vec f$ are bounded and compactly supported.
Note that $H_j\subset\Omega_j$ and
$|H_j|=4^{-(j+1)}$ for all $j\in\mathbb Z_+$. Thus,
\begin{align}\label{direct:input}
\|\vec f\|_{L^p(W)}^p
&= \int_{\mathbb R} |L\vec f|^p
= \int_{\mathbb R} w_0 u^p
= \sum_{j=0}^{J-1} 4^{-(p-1)j} e^{(p-1)\varepsilon j}
\left( \varepsilon  4^j e^{-\frac{\varepsilon}{p'}j} \right)^p |H_j| \notag \\
&= J \frac{\varepsilon^p}{4}
\sim (\delta^{-1}+\varepsilon^{-1}) \varepsilon^p
\sim Q^{-\min\{1,\frac1{p-1}\}}.
\end{align}

Next, we estimate $\|W^{\frac1p}H\vec f\|_{L^{p,\infty}}$.
Since $|LH\vec f|\ge w_1^{\frac1p}|\varphi Hu -H(\varphi u)|$, we have
\begin{equation} \label{est1}
\left\|W^{\frac1p}H\vec f\right\|_{L^{p,\infty}}
= \left\|LH\vec f\right\|_{L^{p,\infty}}
\ge \left\|w_1^{\frac1p}|\varphi Hu -H(\varphi u)|\right\|_{L^{p,\infty}}.
\end{equation}
Note that, for any $I\in\mathcal I_N$,
$\varphi$ is affine on $H_I$.
Thus, for any $x\in H_I$,
\begin{equation}\label{direct:commutator}
\varphi(x)Hu(x)-H(\varphi u)(x)
=\frac1\pi\int_{\mathbb R}
\frac{\varphi(x)-\varphi(y)}{x-y}u(y)\,dy,
\end{equation}
where the integrand is nonnegative
because $\varphi$ is nondecreasing.

For $i\in\{0,\ldots,N-1\}$, let $P_i\in\mathcal I_i$ contain $I$, and let
$S_i$ be the child of $P_i$ that does not contain $I$.
The intervals $S_i$ are pairwise disjoint. For any $y\in S_i$,
$H_{P_i}$ lies between $x$ and $y$, and $|x-y|\le8^{-i}$. Thus,
\begin{equation*}
\frac{\varphi(x)-\varphi(y)}{x-y}
= \frac1{x-y}\int_y^x \varphi'(t)\,dt
\ge 8^i\int_{H_{P_i}}\varphi'(t)\,dt
= \frac1{32} 2^ie^{\frac\varepsilon{p'} i},
\end{equation*}
where the last equality used \eqref{daoshu}.
Note that, for any $j\in\mathbb Z_+$,
if $j\le i$, then $|S_i\cap H_j|=0$, while if $j\ge i+1$,
then $|S_i\cap H_j|=\frac14 8^{-j} 2^{j-i-1}=4^{-j} 2^{-i-3}$. Thus,
\begin{equation*}
\int_{S_i} u
=\varepsilon \sum_{j=i+1}^{J-1} 4^j
e^{-\frac{\varepsilon}{p'}j} |S_i\cap H_j|
= 2^{-i-3} \varepsilon
\sum_{j=i+1}^{J-1}e^{-\frac{\varepsilon}{p'}j}
\sim 2^{-i} e^{-\frac\varepsilon{p'} i}.
\end{equation*}
Restricting the domain of integration on the right-hand side of
\eqref{direct:commutator} to $\bigcup_{i=0}^{N-1} S_i$ yields, for any $x\in H_N$,
\begin{equation*}
\varphi(x)Hu(x)-H(\varphi u)(x)
\gtrsim \sum_{i=0}^{N-1} \left( 2^ie^{\frac\varepsilon{p'} i}\right)
\left( 2^{-i} e^{-\frac\varepsilon{p'} i}\right)
=N.
\end{equation*}
Substituting this estimate into \eqref{est1}, we obtain
$$
\left\|W^{\frac1p}H\vec f\right\|_{L^{p,\infty}}
\gtrsim \left\|w_1^{\frac1p} N\mathbf 1_{H_N}\right\|_{L^{p,\infty}}.
$$
Since $w_1=4^N$ on $H_N$, it follows that
$$
\left\|W^{\frac1p}H\vec f\right\|_{L^{p,\infty}}
\gtrsim 4^{\frac Np} N |H_N|^{\frac1p}
\sim N
\ge Q.
$$
Combining this and \eqref{direct:input}, we obtain
\begin{equation} \label{final}
\frac{\|W^{\frac1p}H\vec f\|_{L^{p,\infty}}}{\|\vec f\|_{L^p(W)}}
\gtrsim Q Q^{\min\{\frac1p,\frac1{p(p-1)}\}}
=Q^{1+\min\{\frac1p,\frac1{p(p-1)}\}}.
\end{equation}

The inequality \eqref{final} readily extends to $m \ge 3$
and can also be expressed in the following uniform form:
for any $t\in[1,\infty)$,
$$
\sup_{[W]_{\mathcal A_p}\le t}
\left\|W^{\frac{1}{p}}H\left(W^{-\frac{1}{p}}\cdot\right)
\right\|_{L^p \to L^{p,\infty}}
\gtrsim t^{1+\min\{\frac1p,\frac1{p(p-1)}\}}.
$$
Both the extension and the reformulation are straightforward;
see \cite[Theorem 3.13]{byyz25}.
This completes the proof of Theorem~\ref{thm:hilbert-model}.
\end{proof}

\noindent\textbf{Acknowledgements}\quad
The authors acknowledge the use of AI tools during the exploratory stage of this project.
All mathematical arguments and proofs in the final manuscript
were checked and written by the authors.

\bigskip
\medskip

\noindent Fan Bu

\medskip

\noindent Department of Mathematics, Faculty of Arts and Sciences,
Beijing Normal University, Zhuhai 519087, The People's Republic of China

\smallskip

\noindent{\it E-mail:} \texttt{fanbu@bnu.edu.cn}

\bigskip

\noindent Dachun Yang and Wen Yuan

\medskip

\noindent Laboratory of Mathematics and Complex Systems (Ministry of Education of China),
School of Mathematical Sciences,  Institute for Advanced Study,
Beijing Normal University, Beijing 100875, The People's Republic of China

\smallskip

\noindent{\it E-mails:} \texttt{dcyang@bnu.edu.cn} (D. Yang)

\noindent\phantom{{\it E-mails:}} \texttt{wenyuan@bnu.edu.cn} (W. Yuan)

\end{document}